\documentclass[10pt,journal]{IEEEtran}

\usepackage{graphicx}          
\usepackage{amsmath}           
\usepackage{amssymb}
\usepackage{algorithm}
\usepackage{algpseudocode}
\usepackage{epstopdf}
\usepackage{cite}
\usepackage{url}
\usepackage{subcaption}
\usepackage[latin1]{inputenc}  

\usepackage{amsthm}

\theoremstyle{definition}

\theoremstyle{remark}
\newtheorem{note}{Remark}
\theoremstyle{plain}

\theoremstyle{definition}
\newtheorem{res}{Result}

\DeclareMathOperator{\diag}{diag}
\DeclareMathOperator{\col}{col}

\DeclareMathOperator{\E}{E}
\DeclareMathOperator{\Cov}{Cov}
\DeclareMathOperator{\Var}{Var}

\DeclareMathOperator*{\argmax}{arg\,max}
\newcommand{\mbf}[1]{\mathbf{#1}}
\newcommand{\mbs}[1]{\boldsymbol{#1}}
\newcommand{\what}[1]{\widehat{#1}}
\newcommand{\wtilde}[1]{\widetilde{#1}}
\newcommand{\wbar}[1]{\overline{#1}}

\newcommand{\0}{\mbf{0}}
\newcommand{\I}{\mbf{I}}
\newcommand{\1}{\mbf{1}}

\newcommand{\T}{\top}

\newcommand{\y}{y}
\newcommand{\ybar}{\wbar{\y}}
\newcommand{\yvec}{\mbf{y}}

\newcommand{\m}{\mu}
\newcommand{\eff}{\delta}
\newcommand{\yvar}{v}
\newcommand{\effvar}{\lambda}
\newcommand{\ycov}{\mbf{V}}

\newcommand{\psub}{\theta}
\newcommand{\pvec}{\mbs{\theta}}

\newcommand{\pder}{\partial}
\newcommand{\score}{\mbs{\phi}}
\newcommand{\g}{\mbf{g}}
\newcommand{\imat}{\mbf{J}}
\newcommand{\rmse}{c}
\newcommand{\conf}{C}

\newcommand{\tva}{\alpha}
\newcommand{\tvb}{\beta}
\newcommand{\tvc}{\gamma}
\newcommand{\snr}{\rho}

\begin{document}

\title{Effect Inference from Two-Group Data \\ with Sampling Bias}
\author{Dave Zachariah and Petre Stoica\thanks{This work has been
    partly supported by the Swedish Research Council (VR) under contract 2018-05040.}}

\maketitle

\begin{abstract}
In many applications, different populations are compared using data
that are sampled in a biased manner. Under sampling biases, 
standard methods that estimate the difference between the population
means yield unreliable inferences. Here we develop an inference method that is resilient to sampling biases and is able to
control the false positive errors under moderate bias levels in
contrast to the standard approach. We demonstrate the method using
synthetic and real biomarker data.
\end{abstract}

\section{Introduction}

In many applications of statistical inference, the aim is
to compare data from different
populations. Specifically, given $n_0$ and $n_1$ samples from two
groups, collected in vectors $\yvec_0$ and $\yvec_1$, the target
quantity is often the difference between their means,
denoted $\eff$, which we call the effect. For instance, in randomized trials and A/B testing,
the data are outcomes from two populations and $\eff$ is the average
causal effect of assigning subjects to a test group `$1$' as compared
to a control group
`$0$'. \cite{Imbens&Rubin2015causal,PearlEtAl2016causal} The standard approach is to use the difference between sample averages in
each group, viz. $\what{\eff} = \wbar{\y}_1 -  \wbar{\y}_0 $, where
$\wbar{\y}_i = \1^\T \yvec_i /n_i$. Confidence intervals
for $\what{\eff}$ can be obtained using Welch's method, which employs an
approximating t-distribution
\cite{Rao1973_linear,Welch1938_significance,Kim&Cohen1998_behrens}. Inferring
$\eff \neq 0$ is equivalent to detecting that the
means of two distributions differ, which is a classical problem in
statistical signal processing \cite{van2013detection,kay1993fundamentals}.

Ideally, the samples from both groups are
representative of their target populations. Then the bias of the estimator,
$$b = \E\left[\what{\eff} \right] - \eff,$$
is zero. However, in nonideal conditions with finite samples this is
not the case, e.g., when some units of the intended populations are
less likely to be included than others. Under such conditions, $b$ decreases with sample sizes $n_0$ and $n_1$ but will
nevertheless be nonzero. Sampling biases increase the risk of
inferring spurious effects when using standard inference methods.

In this paper, we develop an inference method that is resilient to
sampling biases. In contrast to the standard
approach, the proposed method reduces the risk of reporting spurious effect estimates
and is capable of controlling the false positive errors under moderate
biases. The method relies on an effect estimator using a 
fully automatic and data-adaptive regularization. We demonstrate its performance on
both synthetic and real data.

\begin{note}
Code for the method can be found at \url{https://github.com/dzachariah/two-groups-data}
\end{note}

\begin{figure}
  \begin{center}
    \includegraphics[width=0.80\columnwidth]{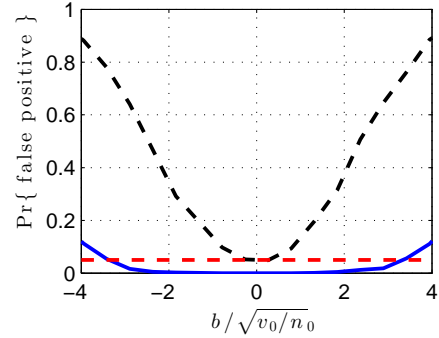}
  \end{center}
  \caption{Probability of false positive error versus bias $b$, when
    $\eff = 0$. Significant
    effects are inferred when the confidence interval excludes the
    zero effect, using Welch's method (black dashed line) and proposed
    method (solid line). Setting $\alpha = 0.05$, the error rate must
    not exceed $5\%$ (red dashed line).  The bias is varied in units of the standard deviation of
    $\ybar_0$ and added to the data from the test group. Data was
    generated using \eqref{eq:model} with $n_0=40$, $n_1=20$, and
    unknown variances $v_0 =  0.3^2$, $v_1 = 0.15^2$ and mean $\m = 1$.}
  \label{fig:estimatorcomparison}
\end{figure}

\section{Problem formulation}



We model the dataset as
\begin{equation} \yvec =
\begin{bmatrix}
\yvec_0 \\
\yvec_1
\end{bmatrix}
\: \sim \: \mathcal{N}\left(
\begin{bmatrix}
\1 & \0 \\
\1 & \1
\end{bmatrix}
\begin{bmatrix}
\m \\
\eff
\end{bmatrix},
\begin{bmatrix}
\yvar_0 \I & \0 \\
\0 & \yvar_1 \I \\
\end{bmatrix}
\right)
\label{eq:model}
\end{equation}
The model based on the Gaussian distribution yields the least favourable distribution for estimating
the unknown effect $\eff$ \cite{stoica2011gaussian}. We model the
effect as a random
variable, where different ranges of values of $\eff$ have different
probabilities. To achieve resilliance to sampling biases, we 
adopt a conservative approach in which nonexistant or
negligible effects are considered to be more probable. Specifically, we
employ the following model:
\begin{equation}
\eff \sim \mathcal{N} (0, \effvar),
\label{eq:modeleff}
\end{equation}
where $\effvar$ is an \emph{unknown} parameter.

Our aim is to derive a confidence interval
$\conf_\alpha(\yvec)$ that contains the unknown $\eff$ with a coverage probability of at least $1-\alpha$. That is,
\begin{equation}
\Pr\Big\{ \: \eff \in \conf_\alpha(\yvec)  \: \Big\}  \; \geq \; 1-\alpha.
\label{eq:coverage}
\end{equation}
The confidence interval is to be centered on an estimator
$\what{\eff}(\yvec)$ and should be resilient to sampling biases. That
is, even if $b\neq 0$ the interval must not indicate 
nonzero effects with a probability greater than
$\alpha$. Fig.~\ref{fig:estimatorcomparison} illustrates the
ability of the method proposed below to ensure \eqref{eq:coverage}
under a range of biases, provided $b$ does not greatly exceed the
dispersion of sample averages, i.e., $\sqrt{v_i/n_i}$.

We will derive a confidence interval using model
\eqref{eq:model} and \eqref{eq:modeleff}, with nuisance
parameters
\begin{equation*}
\pvec = \col\{ \m, \effvar, \yvar_0, \yvar_1 \}.
\end{equation*}

\section{Proposed method}

Let $ \E_\psub[ \eff | \yvec ]$ be the conditional mean of the effect
given the data. Using an estimate $\what{\pvec}$ of the nuisance
parameters, we propose the following effect estimator
\begin{equation}
\begin{split}
\what{\eff}(\yvec) &= \E_\psub[ \eff | \yvec ] \: \Big|_{\psub=\what{\psub}} 
\\
&= \frac{\snr n_1}{\snr n_1 + 1 } (\ybar_1 - \m) \: \Big|_{\psub=\what{\psub}} ,
\end{split}
\label{eq:eb_estimates}
\end{equation}
where we introduce the variable $\snr \equiv \effvar / \yvar_1$ that
can be interpreted as a signal-to-noise ratio, see
\cite{kailath2000linear} for a derivation.

\begin{res}[Cramér-Rao bound]
When the systematic error of $\what{\delta}(\yvec)$ is invariant with respect to
$\pvec$, then the mean-squared error over all possible effects and
data has a Cramér-Rao bound $\E\left[ | \eff - \what{\eff}(\yvec) |^2 \right] \geq \rmse^2_\psub,$
where
\begin{equation}
\rmse^2_\psub = \frac{\snr \yvar_1}{\snr n_1 + 1} +
\frac{\snr^2 n^2_1}{(\snr n_1 + 1)^2} \left(
  \frac{n_0}{\yvar_0} + \frac{n_1}{\yvar_1}\frac{1}{\snr n_1 + 1}
\right)^{-1}.
\label{eq:crb}
\end{equation}
\end{res}
\begin{proof}
See Appendix \ref{app:proofs_crb}.
\end{proof}
\begin{res}[Confidence interval]
Let
\begin{equation}
\conf_\alpha(\yvec) = \big\{ \eff' : |\eff' - \what{\eff}(\yvec)| <
\alpha^{-1/2} \rmse_\psub \big\}.
\label{eq:targetconf}
\end{equation}
When using an efficient
estimator that attains the bound \eqref{eq:crb}, the interval in
\eqref{eq:targetconf} satifies the specified coverage probability \eqref{eq:coverage}.
\end{res}
\begin{proof}
See Appendix \ref{app:proofs_conf}.
\end{proof}

Evaluating $\what{\delta}(\yvec)$ and $\conf_\alpha(\yvec)$ requires
estimates of the nuisance parameters $\pvec$. Here we adopt the maximum
likelihood approach and estimate $\pvec$ using the marginalized data
distribution,
\begin{equation}
p_\psub( \yvec  ) = \int p_\psub( \yvec | \eff )
p_\psub(\eff) d\eff 
\label{eq:marginal}
\end{equation}
It can be shown that \eqref{eq:marginal} is a Gaussian distribution \cite{kailath2000linear}
with mean $\E_\psub[ \yvec  ] = \1 $ and covariance
\begin{equation*}
\Cov_\psub[ \yvec ] = \diag( \yvar_0 \I, \: \effvar \1 \1^\T + \yvar_1
\I ),
\end{equation*}
The estimated parameters are given by
\begin{equation}
\what{\pvec} = \argmax_{\pvec} \: p_\psub(\yvec  ) ,
\label{eq:MMLE}
\end{equation}
which can be shown to yield an asymptotically
efficient estimator \eqref{eq:eb_estimates} \cite[corr.~9]{bar2015bayesian}.

Interestingly, the problem \eqref{eq:MMLE} can be solved by a
one-dimensional numerical search. Begin by defining the variables
\begin{equation*}
\begin{split}
\tva &=  \yvec^\T_1 \yvec_1 
-\m n_1 ( 2 \wbar{\y}_1 -  \m) \\
\tvb &=n^2_1 ( \wbar{\y}_1 - \mu )^2 \\ 
\tvc &= \tva n_1 - \tvb.
\end{split}
\end{equation*}
Note that $\tvc \geq 0$. Then the following result holds.

\begin{res}[Nuisance parameter estimates]
The estimated variances are given by
\begin{equation}
\what{\yvar}_0 = \frac{1}{n_0} \yvec^\T_0 \yvec_0 - \m ( 2 \wbar{\y}_0 -  \m) ,
\label{eq:min_v0}
\end{equation}
\begin{equation}
\begin{split}
\what{\yvar}_1 = \frac{1}{n_1} \frac{\tva + \snr \tvc}{1+\snr n_1},
\end{split}
\label{eq:min_v1}
\end{equation}
which are ensured to be nonnegative, and $\what{\effvar} = \what{\snr} \what{\yvar}_1$, where 
\begin{equation}
\begin{split}
\what{\snr}= \begin{cases} \frac{\tvb - \tva}{\tvc},  & \tvb - \tva \geq 0. \\
0, & \text{otherwise.}
\end{cases}
\end{split}
\label{eq:min_snr}
\end{equation}
All variables in \eqref{eq:min_v0}-\eqref{eq:min_snr} are functions of the mean $\mu$, whose estimate
$\what{\mu}$ is obtained by minimizing the one-dimensional function
\begin{equation}
\begin{split}
f(\m)&= n_0 \ln \what{v}_0 + n_1 \ln (\tva + \what{\snr} \tvc) - (n_1-1) \ln
(1+ \what{\snr} n_1)
\end{split}
\label{eq:finalcost}
\end{equation}
\end{res}
\begin{proof}
See Appendix~\ref{app:proofs_mle}.
\end{proof}

By plugging in $\what{\m}$, $\what{\rho}$, $\what{\yvar}_0$ and
$\what{\yvar}_1$ into \eqref{eq:eb_estimates} and \eqref{eq:targetconf}, we obtain estimates $\what{\eff}(\yvec)$ and $\conf_\alpha(\yvec)$,
respectively. We note that the overall mean $\m$ is fitted to
the data in a nonstandard manner using \eqref{eq:finalcost}, which
yields a fully automatic and data-adaptive regularization of the effect estimator \eqref{eq:eb_estimates}. If the minimizing $\what{\m}$ is such
that $\tvb < \tva$, then the estimated signal-to-noise ratio is
$\what{\snr} = 0$. In this case, the method indicates that the data
is not sufficiently informative to discriminate any systematic
difference from noise.  Consequently, $\what{\eff}(\yvec)$
collapses to zero and $\conf_\alpha(\yvec) = \emptyset$, indicating
a case in which the effect cannot be reliably inferred.

\section{Experimental results}

We demonstrate the proposed inference method using both synthetic and
real data.

\subsection{Synthetic data}

\begin{figure*}
\centering
   \begin{subfigure}[b]{0.42\textwidth}
   \includegraphics[width=1\linewidth]{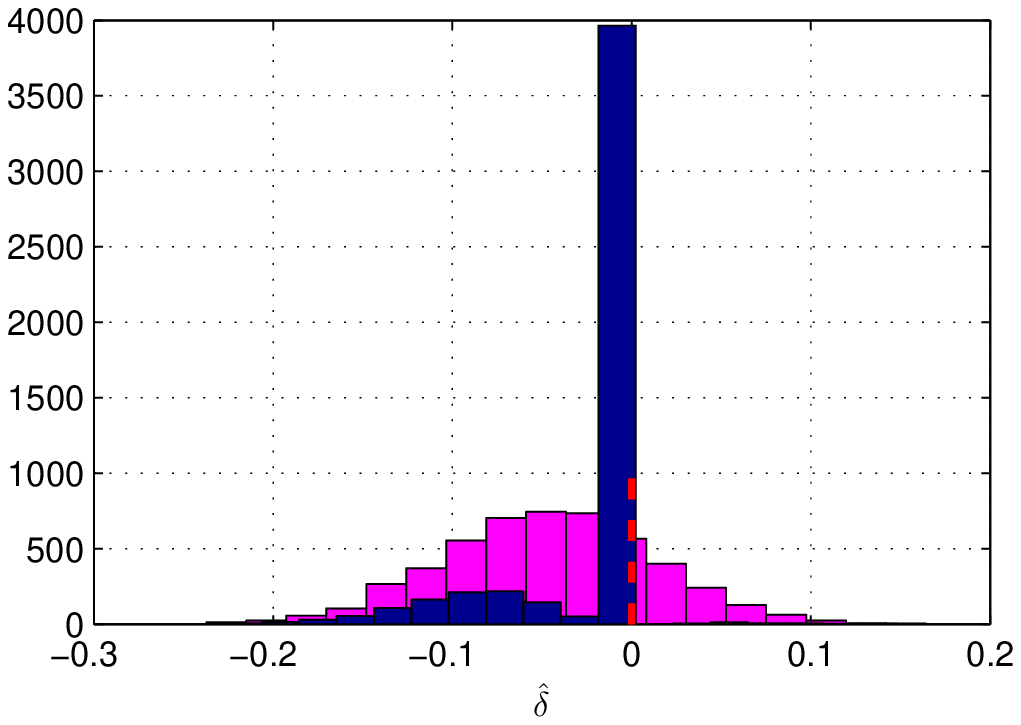}
   \caption{Zero effect $\eff = 0$}
\end{subfigure}
~
\begin{subfigure}[b]{0.42\textwidth}
   \includegraphics[width=1\linewidth]{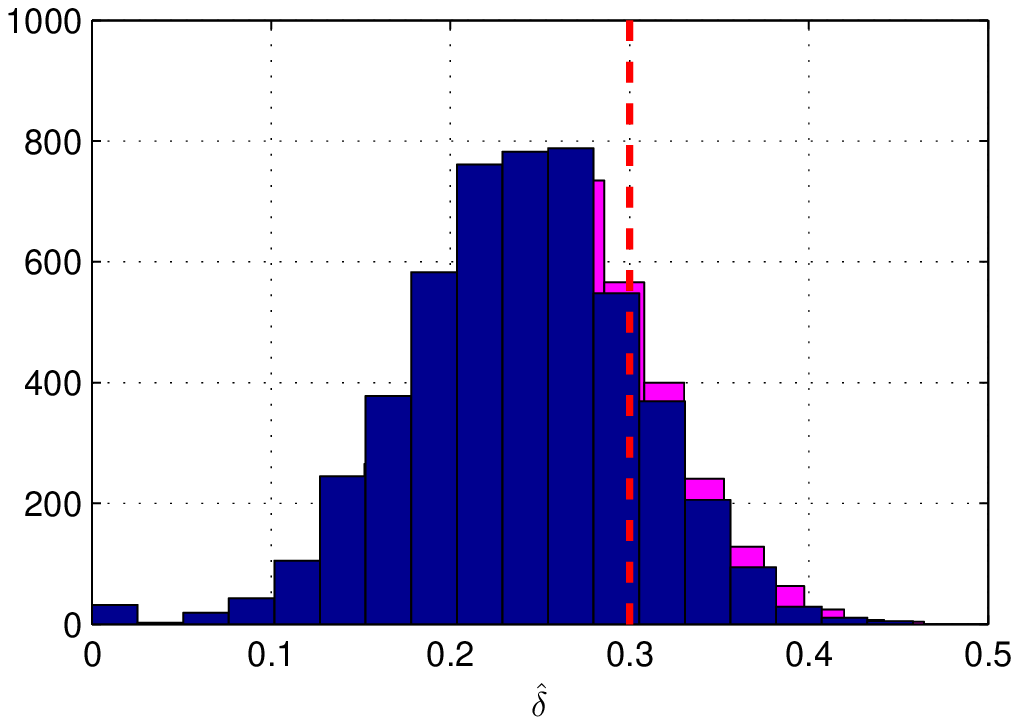}
   \caption{Nonzero effect $\eff = \sqrt{\yvar_0}$}
\end{subfigure}
\caption[TEST]{Distributions of $\what{\eff}$ using standard (pink) and
  proposed (blue) methods under negative 
  bias $b=-\sqrt{\yvar_0/n_0}$. Unknown effect $\eff$ indicated by red
  dashed line. Histograms obtained using $5000$ Monte Carlo realizations.}
\label{fig:histograms}
\end{figure*}

We generate two-group data using the model \eqref{eq:model}
and add a negative bias $b$ to the test group, using the setup
parameters described in Fig.~\ref{fig:estimatorcomparison}. The
adaptive regularization of $\what{\eff}$ is illustrated in Fig.~\ref{fig:histograms}:
when the unknown effect is nonexistent, $\eff=0$, the estimates are
concentrated at zero, despite the bias $b$. As $\eff$ exceeds the
dispersion of the sample averages, however, the regularized and
standard estimators become nearly identical.

\begin{figure*}
\centering
   \begin{subfigure}[b]{0.43\textwidth}
   \includegraphics[width=1\linewidth]{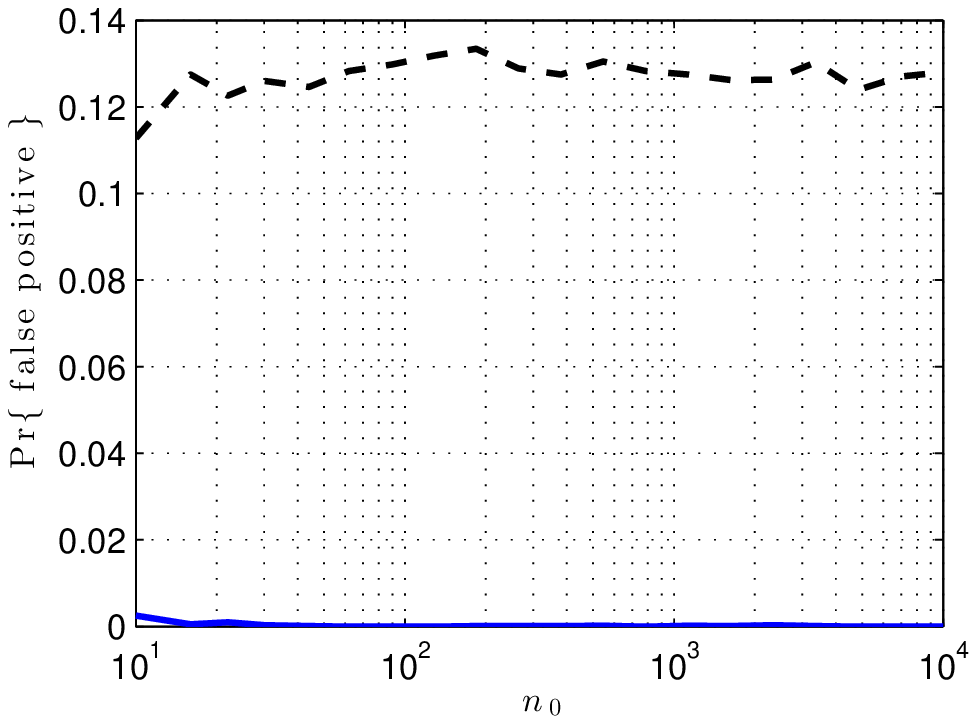}
   \caption{Zero effect $\eff=0$}
\end{subfigure}
~
\begin{subfigure}[b]{0.43\textwidth}
   \includegraphics[width=1\linewidth]{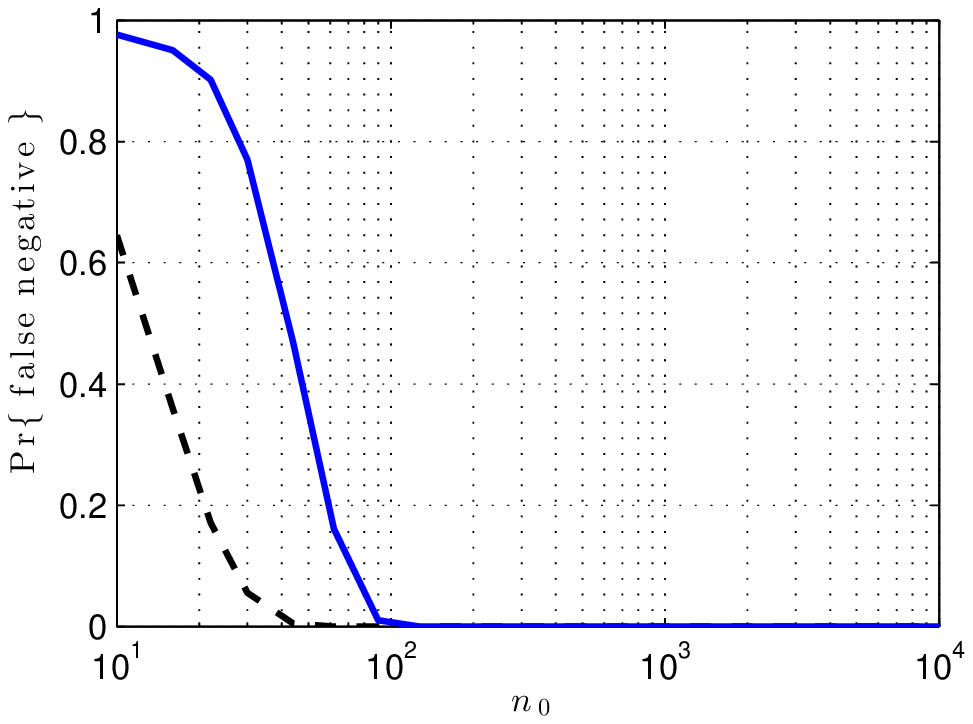}
   \caption{Nonzero effect $\eff = \sqrt{\yvar_0}$}
\end{subfigure}
\caption[TEST]{Probability of false inferences versus  number of
  samples $n_0$, using standard (dashed)
  and proposed (solid) methods. The sample ratio is $n_0/n_1 = 2$ and the 
  bias is $b=-\sqrt{\yvar_0/n_0}$. (a) Probability of false positive
  error, which is targeted to not exceed $\alpha
  = 0.05$. (b) Probability of false negative error, which is the
  complement of the statistical `power'.}
\label{fig:probabilityerror}
\end{figure*}

We report a significant effect estimate when a nonempty interval $\conf_\alpha(\yvec)$ excludes the zero effect. Fig.~\ref{fig:probabilityerror}
illustrates the ability of the proposed method to control the
false positive error probability as $n_0$ increases, in contrast to
the standard method. This is achieved while incurring a 
loss of statistical power that vanishes as the number of samples increases.

\subsection{Prostate cancer data}

\begin{figure*}
  \begin{center}
    \includegraphics[width=1.90\columnwidth]{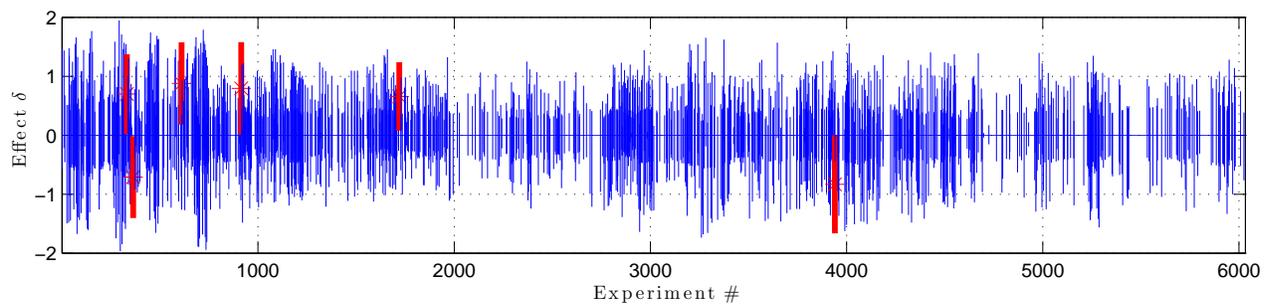}
  \end{center}
  \caption{Confidence intervals $\conf_{\alpha=0.05}(\yvec)$ for 6033
    different experiments using two-group biomarker data ($n_0=50$ and
    $n_1 = 52$). In six
    cases, highlighted in red, the effects were found to be significant
    as the intervals did not contain the zero effect. Note that
    several intervals are empty, indicating cases in which the data is not
    informative enough for the fitted model to discern any systematic effect from the noise.}
  \label{fig:prostate}
\end{figure*}

We now consider real data from $n_0=50$ healthy individuals and
$n_1=52$ individuals with prostate cancer
\cite{Efron2012_largescale,SinghEtAl2002_gene}. The data contains 6033
different biomarker responses. The inferred effects are shown in
Fig.~\ref{fig:prostate}. For 6 markers, the effects were found to be
significant at the $\alpha=0.05$ level. By contrast, the
standard approach using Welch's t-intervals yields 478 genes, but the
inferences are less reliable under sampling biases.


\section{Conclusions}

 We developed a method for inferring effects in two-group data that,
 unlike the standard approach, is resilient to sampling
 biases. The method is able to control the false positive errors under
 moderate bias levels and its performance was demonstrated using both synthetic and real biomarker data.


\appendix

\subsection{The derivation of the Cramér-Rao bound}
\label{app:proofs_crb}

The mean-square error can be decomposed as
\begin{equation}
\begin{split}
\E\left[|\eff - \what{\eff}|^2 \right] &= \E_{y}\left[
  \E_{\eff|y}\left[ |\eff - \wbar{\eff} + \wbar{\eff}  - \what{\eff}
    |^2 \right] \right] \\
&= \E_{y}\left[ \Var[\eff | \yvec] +
  |\wbar{\eff}  - \what{\eff}|^2 \right] \\
&=\frac{\effvar \yvar_1}{\effvar n_1 + \yvar_1} + \E_{y}\left[ |\wbar{\eff}  -
  \what{\eff}|^2 \right].
\end{split}
\end{equation}
where $\wbar{\eff}$ is the conditional mean. Next, define the score
function and the information matrix,
\begin{equation}
\begin{split}
\score \triangleq \pder_\theta \ln p_\psub( \yvec  ) \quad \text{and} \quad \imat = \E_y[\score \score^\T].
\end{split}
\end{equation}
Since the marginal pdf is Gaussian, we can compute $\imat$ using
Slepian-Bangs formula \cite{Stoica&Moses2005_spectral}. It has a block diagonal form
\begin{equation}
\begin{split}
\imat &= \begin{bmatrix}
J_{1,1} & \0 \\
\0    & *
\end{bmatrix},
\end{split}
\end{equation}
where
\begin{equation}
\begin{split}
J_{1,1} &= \1^\T \begin{bmatrix} \yvar^{-1}_0 \I & \0 \\ \0 &
  \ycov^{-1}_1 \end{bmatrix} \1 \\
&=  \yvar^{-1}_0 \1^\T \1 + \1^\T   \ycov^{-1}_1  \1 = \frac{n_0}{\yvar_0} + \frac{n_1}{\effvar n_1 + \yvar_1}
\end{split}
\end{equation}
and $\ycov_1 =\effvar \1 \1^\T + \yvar_1
\I$.

Let $\g \triangleq \E_y[ \score ( \wbar{\eff}- \what{\eff})] $ denote the
correlation between the score function and estimation error. Then we
have the general bound
\begin{equation}
\begin{split}
0 &\leq \E_{y}\left[ | (\wbar{\eff}  - \what{\eff}) - \g^\T \imat^{-1} \score |^2  \right] \\
&= \E_y \left[ | \wbar{\eff}  - \what{\eff}  |^2  \right] - \g^\T \imat^{-1} \g.
\end{split}
\label{eq:covarianceineq}
\end{equation}
In our case, we obtain
\begin{equation*}
\begin{split}
\g &= \int [\pder_\psub \ln p_\psub](  \wbar{\eff} - \what{\eff}) p_\psub d\yvec \\
&= \int \pder_\psub [p_\psub( \wbar{\eff} - \what{\eff})] - p_\psub
[\pder_\psub( \wbar{\eff} - \what{\eff})] d\yvec \\
&=  \pder_\theta [\text{bias}(\pvec)] - \E_{y}[\pder_\theta( \wbar{\eff}  - \what{\eff})] \\
&= - \E_{y}\left[\pder_\theta\left( \frac{\effvar}{\effvar n_1 +
      \yvar_1} \1^\T (\yvec_1 - \mu \1) \right) \right] \\
&= 
-\begin{bmatrix}
 \frac{\effvar}{\effvar n_1 +
      \yvar_1} \1^\T \1 + 0\\
\pder_\effvar \frac{\effvar}{\effvar n_1 +
      \yvar_1} \1^\T \E_y[\yvec_1 - \mu\1] \\
0 \\
\pder_{\yvar_1} \frac{\effvar}{\effvar n_1 +
      \yvar_1} \1^\T \E_y[\yvec_1 - \mu\1] \\
\end{bmatrix} 
= 
-\begin{bmatrix}
 \frac{\effvar n_1}{\effvar n_1 + \yvar_1} \\
0 \\
0 \\
0 \\
\end{bmatrix},
\end{split}
\end{equation*}
where the fourth line follows under the constant bias
assumption. Inserting this expression for $\g$ in \eqref{eq:covarianceineq} yields
\begin{equation}
\begin{split}
\E\left[|\eff - \what{\eff}|^2 \right] \geq \frac{\effvar \yvar_1}{\effvar n_1
  + \yvar_1} + \left(\frac{\effvar n_1}{\effvar n_1 +
    \yvar_1}\right)^2 J^{-1}_{1,1} .
\end{split}
\end{equation}
This completes the proof.

\subsection{The derivation of the confidence interval}
\label{app:proofs_conf}

We have that
\begin{equation}
\begin{split}
\eff \not \in \conf_{\alpha}(\yvec) \quad \Leftrightarrow \quad \alpha
\rmse^{-2}_\psub |\eff - \what{\eff}(\yvec)|^2 \geq 1.
\end{split}
\end{equation}
Let $p(\yvec, \eff) = p_\psub(\yvec| \eff)p_\psub(\eff)$, then
\begin{equation*}
\begin{split}
\Pr\big\{  \eff \not \in \conf_{\alpha}(\yvec) \big\} &= \int_{\eff \not \in \conf_{\alpha}(\yvec)} p(\yvec, \eff) d\eff d\yvec \\
&\leq \int_{\eff \not \in \conf_{\alpha}(\yvec)} \alpha \rmse^{-2}_\psub |\eff - \what{\eff}(\yvec)|^2
p(\yvec, \eff) d\eff d\yvec \\
&\leq \alpha \rmse^{-2}_\psub \E\left[|\eff -
  \what{\eff}(\yvec)|^2\right] = \alpha  \frac{\text{MSE}}{\rmse^2_\psub}.
\end{split}
\end{equation*}
Thus $\Pr\big\{  \delta \in \conf_{\alpha}(\yvec) \big\}  \geq 1- \alpha$ when the
estimator is efficient.



\subsection{The derivation of the concentrated cost}
\label{app:proofs_mle}

Problem \eqref{eq:MMLE} can be formulated equivalently as the
minimization of:
\begin{equation}
\begin{split}
f(\pvec) &= \underbrace{n_0 \ln \yvar_0 +\frac{1}{v_0}\|
\yvec_0 - \m \1 \|^2}_{f_0(\mu, \yvar_0)} +
\underbrace{\ln |\ycov_1| + \| \yvec_1 - \m \1 \|^2_{\ycov^{-1}_1}}_{f_1(\m,\effvar,\yvar_1)}.
\end{split}
\label{eq:costfunction}
\end{equation}
The minimizer
\begin{equation}
\what{\yvar}_0 = \|
\yvec_0 - \m \1 \|^2 / n_0
\end{equation}
is inserted back to yield a concentrated cost function
\begin{equation}
\begin{split}
f_0(\m, \what{\yvar}_0)&= n_0 \ln \what{\yvar}_0 + n_0 
\end{split}
\label{eq:negloglik_0}
\end{equation}
Next, using the Sherman-Morrison and matrix determinant lemmas we can
reparametrize $f_1$ as
\begin{equation}
\begin{split}
f_1(\m, \snr, \yvar)&= \ln (1 + \snr n) + \ln v^n \\
&\quad + \frac{1}{v} \left(
\| \yvec_1 - \mu \1 \|^2 - \frac{\snr |\1^\T (\yvec_1 - \mu \1)
  |^2}{1+\snr n} \right) \\
\end{split}
\label{eq:negloglik_1}
\end{equation}
where we dropped the subindices for notational convenience. 

Using the identities $\tva = \| \yvec_1 - \mu \1 \|^2 $, $ \tvb=|\1^\T
(\yvec_1 - \mu \1) |^2$ and $\tvc = \tva n - \tvb$, the
minimizing $\yvar$ of \eqref{eq:negloglik_1} is found as
\eqref{eq:min_v1}. Inserting the variance estimate back, yields a concentrated cost function
\begin{equation}
\begin{split}
f_1(\m, \snr, \what{\yvar})&= \ln \frac{( \tva + \snr \tvc )^n}{(1+
  \snr n)^{n-1}} + n.
\end{split}
\label{eq:negloglik_1_alt}
\end{equation}
To find the minimizing $\snr \geq 0$, we first consider the stationary point
of 
\begin{equation*}
\begin{split}
\wtilde{f}_1(\m, \snr)&= ( \tva + \snr \tvc )^n(1+
  \snr n)^{-(n-1)}.
\end{split}
\end{equation*}
Taking the derivative with respect to $\snr$, yields the following
condition for a stationary point:
\begin{equation*}
\begin{split}
&n \tvc (\tva + \rho \tvc)^{n-1} (1 + \rho n)^{-n+1} \\
&- (n-1) n (1+\snr
n)^{-n} (\tva + \snr \tvc)^n = 0,
\end{split}
\end{equation*}
or equivalently $\tvc (1 + \snr n) - (n-1)(\tva + \snr \tvc) = 0$.
Solving for $\snr \geq 0$, we obtain the estimate \eqref{eq:min_snr}.

By evaluating the second derivative at this point, we verify that it
is a minimum. 
Inserting \eqref{eq:min_snr} back into \eqref{eq:negloglik_1_alt} and combining with
\eqref{eq:negloglik_0}, we can write \eqref{eq:costfunction} in the
concentrated form \eqref{eq:finalcost}
after omitting irrelevant constants.

\bibliographystyle{ieeetr}
\bibliography{refs_robusteffect}

\end{document}